
\magnification1200
\input amstex.tex
\documentstyle{amsppt}
\nopagenumbers
\hsize=12.5cm
\vsize=18cm
\hoffset=1cm
\voffset=2cm

\footline={\hss{\vbox to 2cm{\vfil\hbox{\rm\folio}}}\hss}

\def\DJ{\leavevmode\setbox0=\hbox{D}\kern0pt\rlap
{\kern.04em\raise.188\ht0\hbox{-}}D}

\def\txt#1{{\textstyle{#1}}}
\baselineskip=13pt
\def\hf{{\textstyle{1\over2}}}
\def\a{\alpha}
\def\d{{\,\roman d}}
\def\e{\varepsilon}\def\E{{\roman e}}

\def\G{\Gamma}

\def\s{\sigma}

\def\={\;=\;}

\def\zt{\zeta(\hf+it)}

\def\E{{\roman e}}

\def\R{\Re{\roman e}\,} 
\def\z{\zeta}

\def\hf{{\textstyle{1\over2}}}
\def\txt#1{{\textstyle{#1}}}

\def\le{\leqslant} \def\ge{\geqslant}
\font\tenmsb=msbm10
\font\sevenmsb=msbm7
\font\fivemsb=msbm5
\newfam\msbfam
\textfont\msbfam=\tenmsb
\scriptfont\msbfam=\sevenmsb
\scriptscriptfont\msbfam=\fivemsb
\def\Bbb#1{{\fam\msbfam #1}}

\def \NN {\Bbb N}
\def \CC {\Bbb C}

\font\ff=cmr8
\def\txt#1{{\textstyle{#1}}}
\baselineskip=13pt

\font\teneufm=eufm10
\font\seveneufm=eufm7
\font\fiveeufm=eufm5
\newfam\eufmfam
\textfont\eufmfam=\teneufm
\scriptfont\eufmfam=\seveneufm
\scriptscriptfont\eufmfam=\fiveeufm
\def\mathfrak#1{{\fam\eufmfam\relax#1}}

\font\tenmsb=msbm10
\font\sevenmsb=msbm7
\font\fivemsb=msbm5
\newfam\msbfam
     \textfont\msbfam=\tenmsb
      \scriptfont\msbfam=\sevenmsb
      \scriptscriptfont\msbfam=\fivemsb
\def\Bbb#1{{\fam\msbfam #1}}

\def \NN {\Bbb N}
\def \CC {\Bbb C}

  \def\rightheadline{{\hfil{\ff
  On a cubic moment of Hardy's function with a shift}\hfil\tenrm\folio}}

  \def\leftheadline{{\tenrm\folio\hfil{\ff
   Aleksandar Ivi\'c }\hfil}}
  \def\emptyheadline{\hfil}
  \headline{\ifnum\pageno=1 \emptyheadline\else
  \ifodd\pageno \rightheadline \else \leftheadline\fi\fi}

\topmatter
\nopagenumbers
\title
ON A CUBIC MOMENT OF HARDY'S FUNCTION WITH A SHIFT
\endtitle
\author   Aleksandar Ivi\'c
 \endauthor

\medskip

\address
Aleksandar Ivi\'c, Katedra Matematike RGF-a
Universiteta u Beogradu, \DJ u\v sina 7, 11000 Beograd, Serbia
\endaddress

\bigskip
\email {
\tt
aleksandar.ivic\@rgf.bg.ac.rs,\;aivic\_2000\@yahoo.com}\endemail
\dedicatory
\enddedicatory
\abstract
{An asymptotic formula for
$$
\int_{T/2}^{T}Z^2(t)Z(t+U)\d t\qquad(0< U = U(T) \le T^{1/2-\e})
$$
is derived, where
$$
Z(t)  := \zt{\bigl(\chi(\hf+it)\bigr)}^{-1/2}\quad(t\in\Bbb R),
\quad \z(s) = \chi(s)\z(1-s)
$$
is Hardy's function. The cubic moment of $Z(t)$ is also discussed,
and a mean value result is presented which supports the author's conjecture that
$$
\int_1^TZ^3(t)\d t \;=\;O_\e(T^{3/4+\e}).
$$
 }
\endabstract
\endtopmatter

\document

 \head
1. Introduction
 \endhead

Let $\z(s)$ denote the Riemann zeta-function, and  as usual
let us define {\it Hardy's function} $Z(t)$ as
$$
Z(t)  := \zt{\bigl(\chi(\hf+it)\bigr)}^{-1/2}\quad(t\in\Bbb R),
\quad \z(s) = \chi(s)\z(1-s).
$$
We have (see e.g., [5, Chapter 1] and [6]), with $s = \s+it$ a complex variable,
$$
\chi(s) = 2^s\pi^{s-1}\sin(\hf \pi s)\G(1-s) =
 {\left(\frac{2\pi}{t}\right)}^{\s+it-1/2}\E^{i(t+\pi/4)}
\left(1+ O\left(\frac{1}{t}\right)\right)\leqno(1.1)
$$
if $t \ge t_0\;(>0)$. In fact, one can obtain a full asymptotic expansion for $\chi(s)$
from the classical Stirling formula for the gamma-function.

\medskip
The function $Z(t)$ (see [8] for an extensive account) is real-valued,
even, smooth, and as $|Z(t)| = |\zt|$, its zeros coincide with
the zeros of $\z(s)$ on the ``critical line'' $\R s = \hf$. It is
indispensable in the study of the zeros of $\z(s)$ on the critical line.

\medskip
The moments of $Z(t)$ are of interest, but since $|Z(t)| = |\zt|$, only odd moments represent
a novelty. M. Korolev [14] and M. Jutila  [10], [11] showed independently that
$$
\int_0^T Z(t)\d t = O(T^{1/4}), \quad \int_0^T Z(t)\d t = \Omega_\pm(T^{1/4}),
$$
thus establishing (up to the value of the numerical constants which are involved) the true
order of the integral in question. This improves on the author's bound $O_\e(T^{1/4+\e})$,
obtained in [7].

\medskip
However, higher odd moments of $Z(t)$, where one also
expects a lot of cancellations due to the oscillatory nature of the function $Z(t)$,
remain a mystery.

\medskip
In [8], equation (11.9), an explicit formula for the cubic
moment of $Z(t)$ was presented. This is
$$
\int\limits_T^{2T}Z^3(t)\d t= 2\pi\sqrt{2\over 3}
\sum_{({T\over 2\pi})^{3/2}\le n\le ({T\over\pi})^{3/2}}
d_3(n)n^{-{1\over6}}\cos\bigl(3\pi n^{2\over 3}+{\txt{1\over8}}\pi\bigr)
+O_\e(T^{3/4+\e}),\leqno(1.2)
$$
where as usual $d_3(n)$ is the divisor function
$$
d_3(n)\= \sum_{k\ell m=n}1\qquad(k,\ell,m,n\in \Bbb N),
$$
generated by $\z^3(s)$ for $\R s >1$.

\medskip
{\bf Remark 1}. Here and later $\e$ denotes positive constants which are arbitrarily
small, but are not necessarily the same ones at each occurrence,
while $a(x) \ll_\e b(x)$ (same as $a(x) = O_\e(b(x)))$ means that
the $|a(x)| \le Cb(x)$ for some $C = C(\e) >0, x\ge x_0$. By $a\asymp b$ we
mean that $a \ll b \ll a$. The symbol $f(x) = \Omega_\pm\bigl(g(x)\bigr)$ means that both
$\limsup\limits_{x\to\infty}f(x)/g(x) > 0$ and $\liminf\limits_{x\to\infty}f(x)/g(x) < 0\,$ holds.

\medskip
The motivation for this investigation is the following problem,
first posed by the author in Oberwolfach 2003 during the conference ``Elementary and
Analytic Number Theory'':
Does there exist  a constant $c$ with $0<c<1$ such that
$$
\int_1^TZ^3(t)\d t \;=\;O(T^c)?\leqno(1.3)
$$
This is equation (11.8) of [8]. So far obtaining (1.3) with any $c < 1$
 is an open problem, but if one considers the
cubic moment of $|Z(t)|$, then it is known that
$$
T(\log T)^{9/4} \;\ll\;  \int_1^T|Z(t)|^3\d t \;=\int_1^T|\zt|^3\d t
\;\ll\; T(\log T)^{9/4},\leqno(1.4)
$$
which establishes the true order of the integral in question. However, obtaining
an asymptotic formula for this integral remains a difficult problem.
The lower bound in (1.4) follows from general results of K. Ramachandra
(see his monograph [15]), and the upper bound is a recent result of S. Bettin, V.
Chandee and M. Radziwi\l {\l} [2].

\medskip
A strong conjecture of the author is that
$$
\int_1^TZ^3(t)\d t \;=\;O_\e(T^{3/4+\e}).\leqno(1.5)
$$
Note that (1.5) would follow from (1.2) and the bound
$$
\sum_{N<n\le N'\le2N}d_3(n)\E^{3\pi in^{2/3}} \ll_\e N^{2/3+\e}.\leqno(1.6)
$$
It may be remarked that the exponential sum in (1.6) is ``pure'' in the sense that the function in
the exponential  does not depend on any parameter as, for example, the sum
$$
\sum_{N<n\le N'\le2N}n^{it}  \;=\;\sum_{N<n\le N'\le2N}\E^{it\log n}
\qquad(1 \le N\ll\sqrt{t}\,),
$$
which appears in the approximation to  $\zt$ (see e.g., Theorem 4.1 of [5]),
depends on the parameter $t$.
However, the difficulty in the estimation of the sum in (1.6) lies in the presence of the
divisor function $d_3(n)$ which, in spite of its simple appearance, is quite difficult
to deal with.

\medskip
{\bf Remark 2}. A natural question is to ask: What is the minimal value of $c$ for which (1.3)
holds ($c = 1+\e$ is trivial)? In other words, to try to find an omega result for the
integral in (1.3).
\medskip
A related and interesting problem  is to investigate
integrals of $Z(t)$ with ``shifts'', i.e., integrals where one (or more) factor $Z(t)$ is
replaced by $Z(t+U)$. The parameter $U$, which does not depend on the variable of
integration $t$, is supposed to be positive and $o(T)$ as $T\to\infty$, where $T$
is the order of the range of integration.

\medskip
Some results on such integrals already exist in the literature. For example,
R.R. Hall [3] proved that, for $U = \a/\log T, \a\ll 1$,
we have uniformly
$$
\eqalign{
\int_0^TZ(t)Z(t+U)\d t &= \frac{\sin\a/2}{\a/2}T\log T + (2\gamma-1-2\pi)T\cos\a/2\cr&
+ O\left(\frac{\a T}{\log T} + T^{1/2}\log T\right).\cr}
\leqno(1.7)
$$
In (1.7) $\gamma= -\G'(1) = 0.5772157\ldots\,$ is Euler's constant.
M. Jutila [12] obtained recently an asymptotic formula for the left-hand side of (1.5) when
$0 < U \ll T^{1/2}$.

\medskip
S. Bettin [1] evaluated asymptotically a related integral, namely
$$
\int_0^T \z(1/2+U+it)\z(\hf-V-it)\d t
$$
under the condition that
$$
U = U(T) \in \CC,\;V = V(T) \in \CC,\; \R U \ll 1/\log T, \; \R V \ll 1/\log T.
$$

\medskip
S. Shimomura [16] dealt with the quartic moment
$$
\int_0^T Z^2(t)Z^2(t+U)\d t,\leqno(1.8)
$$
under certain conditions on the real parameter $U$, such that $|U| + 1/\log T\to0$ as
$T\to\infty$. When $U = 0$, Shimomura's expression for (1.8) reduces to
$$
\int_0^T|\zt|^4\d t = \int_0^T Z^4(t)\d t = \frac{1}{2\pi^2}T\log^4T + O(T\log^3T).
\leqno(1.9)
$$
The (weak) asymptotic formula (1.9) is a classical result of A.E. Ingham [4] of 1928.
For the sharpest known asymptotic formula for the integral in (1.9), see
the paper of Ivi\'c--Motohashi [9].

\medskip
\head
2. Statement of results
\endhead
\medskip

\medskip
The first purpose of this work is to investigate another integral with the shift of $Z(t)$,
namely the cubic moment
$$
\int_{T/2}^T Z^2(t)Z(t+U)\d t\leqno(2.1)
$$
for a certain range of the ``shift parameter'' $U = U(T) >0$.
We shall prove the following

\medskip
THEOREM 1. {\it For $0 < U = U(T) \le T^{1/2-\e}$ we have, uniformly in $U$},
$$
\eqalign{&
\int_{T/2}^{T}Z^2(t)Z(t+U)\d t = O_\e(T^{3/4+\e}) \,+
\cr& \,+
2\pi\sqrt{2\over 3}\sum_{T_1\le n\le T_0}h(n,U)n^{-1/6+iU/3}
\exp(-3\pi in^{2/3}-\pi i/8)\Bigl\{1 + K(n,U)\Bigr\}.
\cr}\leqno(2.2)
$$
{\it Here} ($d(n)$ {\it is the number of divisors of $n$})
$$
h(n,U) := n^{-iU}\sum_{\delta|n}d(\delta)\delta^{iU},
\; T_0 := {T^{3/2}\over\sqrt{8\pi^3}}, \;T_1 := {{(\frac{T}{2})}^{3/2}\over\sqrt{8\pi^3}},
\leqno(2.3)
$$
$$
K(n,U) := d_2U^2n^{-2/3} + \cdots + d_kU^kn^{-2k/3}
+ O_k(U^{k+1}n^{-2(k+1)/3})\leqno(2.4)
$$
{\it for any given integer $k\ge 2$, with effectively computable constants
$d_2, d_3,\ldots\;$.}

\medskip
{\bf Remark 3}. When $U\to0+$, the expressions in (2.2)--(2.4) tend to (1.2),
since the integral is real and $\lim\limits_{U\to0+}h(n,U) = d_3(n)$. Thus
Theorem 1 is a generalization of (1.2) for a relatively wide range of $U$.

\medskip
{\bf Remark 4}. By methods similar to those used in proving Theorem 1 one can deal also
with the ``conjugate'' integral
$$
\int_{T/2}^{T}Z^2(t+U)Z(t)\d t\qquad(0 < U \le T^{1/2-\e}).
$$

\medskip
{\bf Remark 5}. Analogously to (1.5) one may conjecture that, uniformly for
$0 < U \le T^{1/2-\e}$,
$$
\int_0^T Z^2(t)Z(t+U)\d t \;=\;O_\e(T^{3/4+\e}).
$$

\medskip
The second aim of this work is to demonstrate, although (1.5) and (1.6) seem at present intractable,
that it is possible to show that the conjecture (1.6) holds in a
certain mean value sense. To this end let
$$
S(\a,N) := \sum_{N<n\le N'\le2N}d_3(n)\E^{\a in^{2/3}}\qquad(\a\in\Bbb R),\leqno(2.5)
$$
so that the exponential sum in (1.6) is $S(3\pi,N)$.
Then we have

\medskip
THEOREM 2. {\it Every finite interval $\,[A,B]\;(0 < A < B)\,$ contains at least
one point $C$ such that}
$$
S(C,N) \;\ll_\e\; N^{2/3}\log^{9/2}N.\leqno(2.6)
$$
\medskip
\head
3. Proof of the Theorems
\endhead
\medskip
We begin with the proof of Theorem 1.
Suppose that $0 < U \le T^{1/2-\e}$ and let
$$
I := \int\limits_{T/2}^T Z^2(t)Z(t+U)\d t
= \frac{1}{i}\int\limits_{1/2+iT/2}^{1/2+iT}\z^2(s)\z(s+iU){(\chi^2(s)\chi(s+iU))}^{-1/2}\d s.
\leqno(3.1)
$$
The procedure of writing a real-valued integral like a complex integral is fairly
standard in analytic number theory. For example,  see
 the proof of Theorem 7.4 in E.C. Titchmarsh's monograph [17]
on $\z(s)$ and M. Jutila's work [12]. It allows one flexibility
by suitably deforming the contour of integration
in the complex plane. Incidentally, this method of proof is different from the
proof of (1.2) in [5], which is based on the use of approximate functional
equations.

\medskip
Let
$$
h(n,U) := n^{-iU}\sum_{\delta|n}d(\delta)\delta^{iU},
\quad N_0 = N_0(T,U) \;:= \;\sqrt{\frac{T^2(T+U)}{8\pi^3}},\leqno(3.2)
$$
so that $h(n,U) \ll_\e n^\e$. We have
$$
I \= I_1 + I_2,
$$
say, where
$$
\eqalign{
I_1 &:= \frac{1}{i}\int_{1/2+iT/2}^{1/2+iT}\,\sum_{n\le N_0}h(n,U)n^{-s}
{(\chi^2(s)\chi(s+iU))}^{-1/2}\d s,\cr
I_2&:= \frac{1}{i}\int_{1/2+iT/2}^{1/2+iT}\left\{\z^2(s)\z(s+iU)
- \sum_{n\le N_0}h(n,U)n^{-s}\right\}{(\chi^2(s)\chi(s+iU))}^{-1/2}\d s,\cr}
\leqno(3.3)
$$
and we shall keep in mind that for $\R s>1$, since $d(n)$ is generated by $\z^2(s)$, 
$$
\sum_{n=1}^\infty h(n,U)n^{-s} = \sum_{k=1}^\infty\sum_{m=1}^\infty
d(k)k^{iU}(km)^{-s}(km)^{-iU}
= \z^2(s)\z(s+iU).
\leqno(3.4)
$$
To be able to exploit (3.4) let henceforth $c = 1+\e$. Then, by Cauchy's theorem,
$$
I_2 = \frac{1}{i}\left\{\int_{1/2+iT/2}^{c+iT/2} + \int_{c+iT}^{1/2+iT}+
\int_{c+iT/2}^{c+iT}\right\} = \frac{1}{i}(J_1+J_2+J_3),
$$
say. For $T/2 \le t \le T, \,\hf \le \s \le c, \,s = \s + it, \,0 < U \ll T^{1/2-\e}\,$ we have
$$
\chi^2(s)\chi(s+iU) = \E^{\frac{3\pi i}{4}}\left({2\pi\over t}\right)^{2\s+2it-1}
\left({2\pi\over t+U}\right)^{\s+it+iU-\frac{1}{2}}\E^{3it+iU}\left\{1+ O\left({1\over t}\right)\right\},
\leqno(3.5)
$$
so that uniformly in $U$
$$
\chi^2(s)\chi(s+iU) \;\asymp\;T^{3/2-3\s}\qquad(t \asymp T).
$$
This follows from the asymptotic formula (1.1), and actually
the term $O(1/t)$ can be replaced by an asymptotic expansion. We also have the
convexity bound (see e.g., Chapter 1 of [5])
$$
\z(s) \ll_\e |t|^{(1-\s)/2+\e}\qquad(\hf \le \s \le c).\leqno(3.6)
$$
Stronger bounds than (3.6) are available, but this bound is sufficient for our
present purposes.

\medskip
By using (3.5), (3.6) and the trivial bound
$$
\sum_{n\le N_0}h(n,U)n^{-s} \ll_\e N_0^{1-\s+\e} \ll_\e T^{(3-3\s)/2+\e}
\quad(\hf \le \s \le c),
$$
it follows readily that
$$
J_1 + J_2 \ll_\e \int_{1/2}^{c}T^{(3-3\s)/2+\e}T^{(6\s-3)/4}\d\s \ll_\e T^{3/4+\e}.
\leqno(3.7)
$$
In $J_3$ we have $s = c +it,\, T/2\le t \le T$. Hence
$$
\eqalign{
J_3 &= i\sum_{n>N_0}h(n,U)n^{-c}\int_{T/2}^T n^{-it}\E^{-3\pi i/8}\left({t\over2\pi}\right)^{c+it-1/2}
\cr&
\times \left({t+U\over2\pi}\right)^{(2c+2it+2iU-1)/4}\E^{-(3it+iU)/2}\d t\cr&
+ O\left(\sum_{n>N_0}|h(n,U)|n^{-c}T^{(6c-3)/4}\right).\cr}
$$
The expression in the $O$-term is clearly $\ll_\e T^{3/4+\e}$, and it remains to estimate
$$
\sum_{n>N_0}h(n,U)n^{-c}\E^{-iU/2}\int_{T/2}^T K(t)\E^{if_n(t)}\d t,\leqno(3.8)
$$
with
$$
K(t) := \left({t\over2\pi}\right)^{c-1/2}\left({t+U\over2\pi}\right)^{(2c-1)/4}
\ll_\e T^{3/4+\e}
$$
for $c = 1+\e$. In (3.8) we have
$$
\eqalign{
f_n(t):&= t\log{t\over2\pi} + \hf (t+U)\log{t+U\over2\pi}- \txt{3\over2}t - t\log n,\cr
f_n'(t) &=  \log{t\over2\pi} + \hf \log{t+U\over2\pi}-\log n,\cr
f_n''(t) &= {1\over t} + {1\over2(t+U)} = {3t+2U\over 2t(t+U)} \asymp {1\over T}.
\cr}
\leqno(3.9)
$$
In the sum in (3.8) set $n = [N_0] + \nu, \nu \in\NN$. By the second derivative test
(Lemma 2.2 of [5]) the contribution of the term $\nu=1$ is $\ll T^{\e-1/4}$.
The terms $\nu\ge 2$ are estimated by the use of the first derivative test
(Lemma 2.1 of [5]). For $\nu \gg N_0$ it is found that $|f_n'(t)| \gg 1$, while
for $\nu \ll N_0$ one has, for a suitable $c>0$,
$$|f_n'(t)| \;=\; \log\frac{[N_0]+\nu}{\frac{t}{2\pi}\sqrt{\frac{t+U}{2\pi}}}
\;\ge\; \log\left(1 +c\nu T^{-3/2}\right)\gg \nu T^{-3/2}.
$$
Thus the contribution of the terms with $\nu \ge 2$ will be
$$
\eqalign{
&
\ll_\e T^{3/4+\e}\left(\sum_{2\le\nu\ll N_0}+ \sum_{\nu\gg N_0}\right)|h([N_0]+\nu,U)|
n^{-c}\max_{T/2\le t \le T}{|f_n'(t)|}^{-1}\cr&
\ll_\e T^{3/4+\e}\left(\sum_{2\le\nu\ll N_0}(N_0+\nu)^{-1-\e}T^{3/2}\nu^{-1}
+ \sum_{\nu\gg N_0}\nu^{\e/2-1-\e}\right)\cr&
\ll_\e T^{3/4+\e},\cr}
$$
since $N_0 \asymp T^{3/2}$.

\medskip
It  remains to deal with $I_1$. We have
$$
\eqalign{
I_1 & = \E^{-3\pi i/8}\int_{T/2}^T\sum_{n\le N_0}h(n,U)n^{-1/2}\E^{-iU/2}\E^{if_n(t)}\d t
+ O_\e(T^{3/4+\e})\cr&
= \E^{-iU/2}\E^{-3\pi i/8}\sum_{n\le N_0}h(n,U)n^{-1/2}{\Cal J}_n(T) + O_\e(T^{3/4+\e}),\cr}
\leqno(3.10)
$$
say, where
$$
{\Cal J}_n(T) \;:=\; \int_{T/2}^T\E^{if_n(t)}\d t.\leqno(3.11)
$$
The exponential integral in (3.11) has a saddle point $t_n$, the root
of the equation $f_n'(t) = 0$, which occurs for
$$
t_n^2(t_n+U) = 8\pi^3 n^2.\leqno(3.12)
$$
It is immediately seen that $t_n \asymp n^{2/3} \asymp T$, and that the equation
(3.12) is cubic in $t_n$ for a given $n$ and $U$, and thus it is not easily solvable.
However, its solution can be found asymptotically by an iterative process.
First from (3.12) we have
$$
t_n^3 = 8\pi^3 n^2 + O(UT^2),
$$
which yields
$$
t_n = 2\pi n^{2/3} + O(U).\leqno(3.13)
$$
As we require $T/2 \le t_n\le T$, from (3.12) we infer that $N_1 \le n \le N_0$,
where $N_0$ is defined by (3.2) (this is the reason for the choice of $N_0$) and
$$
N_1 := \sqrt{{(T/2)^2(T/2+U)\over 8\pi^3}}.
$$
Setting $n = [N_1]-\mu, 1 \le \mu < N_1$ and proceeding similarly as in the estimation
of the sum with $\nu\ge2$ after (3.9), we obtain that the contribution of
$n\le N_1$ to $I_1$ in (3.10) will be $O_\e(T^{3/4+\e})$. Since ${\Cal J}_n(T)
\ll T^{1/2}$ (see (3.16) and (3.17))
it transpires that the condition
$N_1 \le n \le N_0$ may be replaced by
$T_1\le n \le T_0$, where $T_0, T_1$ are defined by (2.3). In this process the total
error term will be $O_\e(T^{3/4+\e})$. Therefore,
for the range $0 < U \le T^{1/2-\e}$, we have uniformly
$$
I_1 = \E^{-iU/2}\E^{-3\pi i/8}\sum_{T_1\le n\le T_0}h(n,U)n^{-1/2}{\Cal J}_n(T) +
 O_\e(T^{3/4+\e}),\leqno(3.14)
$$
and the range of summation in (3.14) is easier to deal with than the range of
summation $N_1 <  n \le N_0$. At this point we shall evaluate ${\Cal J}_n(T)$ by
one of the theorems related to exponential integrals with a saddle point,
e.g., by the one on p. 71 of the monograph by Karatsuba--Voronin [13]. This says that,
if $f(x) \in C^{(4)}[a,b]$,
$$
\eqalign{
\int_a^b \E^{2\pi if(x)}\d x &=
\E^{\pi i/4}{\E^{2\pi if(c)}\over \sqrt{f''(c)}} + O(AV^{-1})\cr&
+\,O\Bigl(\min({|f'(a)|}^{-1},\sqrt{a}\,)\Bigr) + O\left(\min({|f'(b)|}^{-1},\sqrt{b}\,)\right),
\cr}
\leqno(3.15)
$$
where
$$
\eqalign{
&0 < b - a \le V,\; f'(c) =0, \;a\le c \le b,\cr&
f''(x) \asymp A^{-1},\; f^{(3)}(x) \ll (AV)^{-1},\; f^{(4)}(x) \ll A^{-1}V^{-2}.
\cr}
$$
We shall apply (3.15) with
$$
f(x) = {f_n(x)\over2\pi},\;c = t_n,\; a = \hf T, \;b=T,\; V = \hf T, f''(x) \asymp {1\over T},
$$
so that $A = T, f^{(3)}(x) \ll T^{-2}, f^{(4)}(x) \ll T^{-3}$, which is needed. The term
$O(AV^{-1})$ in (3.15) will make a contribution which is $O_\e(T^{3/4+\e})$, and the
same assertion will hold for the other two error terms in (3.15). This will lead to
$$
I_1 = \E^{-iU/2}\E^{-3\pi i/8}\sum_{T_1\le n\le T_0}h(n,U)n^{-1/2}\E^{\pi i/4}
{\E^{if_n(t_n)}\over \sqrt{{f_n''(t_n)\over2\pi}}} + O_\e(T^{3/4+\e}).
\leqno(3.16)
$$
Since, by (3.9) and (3.12),
$$
\eqalign{
\E^{if_n(t)} &= \exp\Bigl\{i\Bigl(t_n\log{t_n\over2\pi} + \hf(t_n+U)\log{(t_n+U)\over2\pi}
-\txt{3\over2}t_n-t_n\log n\Bigr)\Bigr\}\cr&
=\exp\Bigl\{\hf it_n\log{t_n^2(t_n+U)\over8\pi^3}+\hf Ui\log{(t_n+U)\over2\pi}
-\txt{3\over2}it_n - it_n\log n\Bigr\}\cr&
=\exp\Bigl\{\hf iU\log{(t_n+U)\over2\pi}- \txt{3\over2}it_n\Bigr\}\cr&
= \E^{-3it_n/2}{\left({t_n+U\over2\pi}\right)}^{iU/2},
\cr}
$$
it follows that
$$
\eqalign{&
I_1 = O_\e(T^{3/4+\e})\cr&
+ \sqrt{2\pi}\E^{-iU/2}\E^{-\pi i/8}\sum_{T_1\le n\le T_0}h(n,U)n^{-1/2}
\sqrt{{2t_n(t_n+U)\over3t_n+2U}}\E^{-3it_n/2}{\left({t_n+U\over2\pi}\right)}^{iU/2}.\cr}
\leqno(3.17)
$$
In view of (3.12), (3.13) and $n\asymp T^{3/2}$ we have
$$
\eqalign{&
n^{-1/2}\sqrt{{2t_n(t_n+U)\over3t_n+2U}} =
n^{1/2}\sqrt{{2(2\pi)^3\over t_n(3t_n+2U)}}\cr&
= \sqrt{2(2\pi)^3)}n^{1/2}t_n^{-1}\Bigl(3^{-1/2} + O(U/T)\Bigr)\cr&
= \sqrt{{2\over3}(2\pi)^3}n^{1/2}\Bigl((8\pi^3n^2)^{1/3} + O(U^{-1})\Bigr)^{-1} + O(UT^{-5/4})
\cr&
= 2\sqrt{{\pi\over3}}n^{-1/6} + O(UT^{-5/4}).\cr}
$$
It remains to evaluate
$$
{\Cal A} := {\Cal A}(n,U) = \E^{-iU/2}\E^{-3it_n/2}{\left({t_n+U\over2\pi}\right)}^{iU/2}.
\leqno(3.18)
$$
To achieve this we need a more precise expression for $t_n$ than (3.13).
Putting (3.13) in (3.12) we first have
$$
\eqalign{
t_n^3 &= 8\pi^3n^2 - Ut_n^2 = 8\pi^3n^2 - 4\pi^2Un^{4/3} + O(n^{2/3}U^2),\cr
t_n &= {(8\pi^3n^2 - 4\pi^2Un^{4/3})}^{1/3}\left(1 + O(U^2n^{-4/3})\right)\cr&
= 2\pi n^{2/3} - \txt{1\over3}U + O(U^2n^{-2/3}).\cr}
$$
This yields the second approximation
$$
t_n = 2\pi n^{2/3} - \txt{1\over3}U + O(U^2n^{-2/3}).\leqno(3.19)
$$
Inserting (3.19) in (3.12) we shall obtain next
$$
t_n =2\pi n^{2/3} - \txt{1\over3}U - {U^2n^{-2/3}\over18\pi} + O(U^3n^{-4/3}).\leqno(3.20)
$$
In fact, (3.20) is just the beginning of an asymptotic expansion for $t_n$ in which
each term is by a factor of $Un^{-2/3}$ of a lower order of magnitude than the
preceding one. This can be established by mathematical induction.
Then we see that ($c$'s and $d$'s are effectively computable constants)
$$
\eqalign{
\log{\Cal A}&= i\Bigl\{-\hf U - {3\over2}(2\pi)n^{2/3} +\hf U + {U^2n^{-2/3}\over12\pi}
+ c_3U^3n^{-4/3} + O(U^4n^{-6/3})\cr&
+ \hf U\log\Bigl(n^{2/3}+{U\over3\pi} - {U^2n^{-2/3}\over36\pi^2} + O(U^3n^{-4/3})\Bigr)\Bigr\}
\cr&
= -3\pi in^{2/3} + \txt{1\over3}i U\log n + ic_2U^2n^{-2/3} + ic_3U^3n^{-4/3}
+ O(U^4n^{-6/3}),\cr&\cr}
$$
and therefore
$$
{\Cal A}= n^{iU/3}\E^{-3\pi in^{2/3}}\left\{1 + id_2U^2n^{-2/3} + id_3U^3n^{-4/3}
+ O(U^4n^{-6/3})\right\}, \leqno(3.21)
$$
and we remark that the last term in curly brackets admits an asymptotic expansion. Since
$U \ll T^{1/2-\e}$ by assumption, $U^2n^{-2/3} \ll T^{-2\e}$ and thus each term in the
expansion will be by a factor of at least $T^{-2\e}$ smaller than the preceding one. One sees
now why the interval of integration in (2.1) and (2.2) is $[T/2,\,T]$ and
not the more natural $[0,\,T]$,  because the terms $d_jU^jn^{-2j/3}$
become large if $n$ is small, and $K(n,U)$ is then also large.

Inserting (3.21) (as a full asymptotic expansion) and (3.18) in (3.17),
we obtain the assertion of Theorem 1, having in mind
that we have shown that 
$$
|I_2| + |I_3| \;\ll_\e\; T^{3/4+\e}.
$$

\medskip
To prove Theorem 2 note that, with $S(\a,N)$ defined by (2.5), we have
$$
\eqalign{
&
\int_A^B|S(\a,N)|^2\d\a = \int\limits_A^B\left(\sum_{n\asymp N}d_3^2(n) + \sum_{m\ne n\asymp N}
d_3(m)d_3(n)\E^{i\a(m^{2/3}-n^{2/3})}\right)\d\a
\cr&
= (B-A)\sum_{n\asymp N}d_3^2(n) + O\left(\sum_{m\ne n\asymp N}
\frac{d_3(m)d_3(n)}{|m^{2/3}-n^{2/3}|}\right).
\cr}\leqno(3.22)
$$

Note that the function $d_3(n)$ is multiplicative and $d_3(p) = 3$. Thus
for $\R s > 1$ and $p$ a generic prime we have
$$
\eqalign{&
\sum_{n=1}^\infty d_3^2(n)n^{-s} = \prod_p\Bigl(1 + d_3^2(p)p^{-s} + d_3^2(p^2)p^{-2s}+\cdots\Bigr)
\cr& =  \z^9(s)\prod_p(1-p^{-s})^9\prod_p\Bigl(1 + 9p^{-s} + d_3^2(p^2)p^{-2s}+\cdots\Bigr)
\cr&
= \z^9(s)F(s),
\cr}\leqno(3.23)
$$
where $F(s)$ is a Dirichlet series that is absolutely convergent for $\R s > 1/2$.
Thus the series in (3.23) is dominated by $\z^9(s)$ and it follows, by a simple
convolution argument, that
$$
\sum_{n\le x}d_3^2(n)\; =\; C_1x\log^8x \;+\; O(x\log^7x)\leqno(3.24)
$$
for a constant $C_1>0$.

\medskip
Next, by using the trivial inequality
$$
|ab| \;\le \;\hf(a^2 + b^2)
$$
and (3.24),
it is seen that the double sum in (3.22) in the $O$-term is, by symmetry,
$$
\eqalign{
& \ll N^{1/3}\sum_{m\ne n\asymp N}\frac{d^2_3(m)+d^2_3(n)}{|m-n|}
\cr&
=2N^{1/3}\sum_{m\asymp N}d_3^2(m)\sum_{n\ne m,n\asymp N}\frac{1}{|n-m|}
\cr&
\ll N^{4/3}\log^9N,\cr}
$$
since, for a fixed $m$ such that $m \asymp N$,
$$
\sum_{n\ne m,n\asymp N}\frac{1}{|n-m|} \;\ll\; \log N.
$$
From this estimate and (3.24)  it follows that
$$
\int_A^B|S(\a,N)|^2\d\a \;\ll\; N^{4/3}\log^9N,\leqno(3.25)
$$
which implies the assertion of Theorem 2.
A further discussion on the cubic moment of
$Z(t)$ is to be found in Chapter 11 of the author's monograph [8].

\medskip
{\bf Remark 6}. An interesting problem is to determine the true order of the integral
in (3.25).

\medskip
{\bf Remark 7}. If the divisor function $d_3(n)$ is not present in the definition (2.5),
namely if one considers the sum
$$
T(\a, N) := \sum_{N<n\le N'\le 2N}\E^{\a in^{2/3}} \qquad(\a\ne 0),
$$
then the problem is much easier. By using Lemma 1.2 and Lemma 2.1 of [5] it easily follows that
$$
T(\a, N) \;\ll |\a|^{-1}N^{1/3}.
$$

\bigskip
\vfill
\eject
\topglue1cm
\medskip
\Refs
\medskip

\item{[1]} S. Bettin,
{\it The second moment of the Riemann zeta function with unbounded shifts},
 Int. J. Number Theory {\bf6}(2010), no. 8, 1933-1944.
\item{[2]} S. Bettin, V.
Chandee and M. Radziwi\l \l, {\it The mean square of the product of the Riemann
zeta function with Dirichlet polynomials}, to appear in
Journal f\"ur die reine und angewandte Mathematik, available online at DOI: 10.1515/crelle-2014-0133.

\item{[3]} R.R. Hall, {\it A new unconditional result about large spaces
between zeta zeros}, Mathematika {\bf53}(2005), 101-113.

\item{[4]} A.E. Ingham, {\it Mean-value theorems in the
theory of the Riemann zeta-function},
Proc. Lond. Math. Soc. (2){\bf27}(1928), 273-300.

\item{[5]} A. Ivi\'c, {\it The Riemann zeta-function}, John Wiley \&
Sons, New York, 1985 (2nd ed. Dover, Mineola, New York, 2003).

\item {[6]} A. Ivi\'c,  {\it Mean values of the Riemann zeta-function},
LN's {\bf 82},  Tata Inst. of Fundamental Research,
Bombay,  1991 (distr. by Springer Verlag, Berlin etc.).

\item {[7]} A. Ivi\'c, {\it On the integral of Hardy's function},  Arch. Mathematik
{\bf83}(2004), 41-47.

\item{[8]} A. Ivi\'c, {\it The theory of Hardy's $Z$-function},
Cambridge University Press, Cambridge, 2012, 245pp.

\item{[9]} A. Ivi\'c and Y. Motohashi, {\it On the fourth power moment of the
    Riemann zeta-function},
    Journal of Number Theory {\bf51}(1995), 16-45.

\item{[10]} M. Jutila, {\it Atkinson's formula for Hardy's function}, J. Number Theory
{\bf129}(2009), no. 11, 2853-2878.

\item{[11]} M. Jutila, {\it An asymptotic formula for the primitive of Hardy's
function},   Arkiv Mat.  {\bf49}(2011), 97-107.

\item{[12]} M. Jutila, {\it The mean value of Hardy's function in short intervals},
Indagationes Math. 2015 (special vol. dedicated to 125th anniversary of J.G. van
der Corput),  {\tt http://dx.doi.org/10.\break 1016/j.indag.2015.08.006}

\item{[13]} A.A. Karatsuba and S.M. Voronin, {\it The Riemann
zeta-function}, Walter de Gruyter, Berlin--New York, 1992.

\item{[14]}  M.A. Korolev, {\it On the integral of Hardy's function $ Z(t)$},
 Izv. Math. {\bf72}, No. 3, (2008), 429-478; translation from  Izv.
Ross. Akad. Nauk, Ser. Mat. {\bf72}, No. 3, (2008), 19-68.

\item{[15]} K.  Ramachandra, {\it On the mean-value and omega-theorems
for the Riemann zeta-function}, LN's {\bf85}, Tata Inst. of Fundamental Research
 Bombay, 1995 (distr. by Springer Verlag, Berlin etc.).

\item{[16]} S. Shimomura, {\it Fourth moment of the Riemann Zeta-function with a shift
along the real line}, Tokyo J. Math. {\bf36}(2013), 355-377.

\item{[17]}  E.C. Titchmarsh, {\it The theory of the Riemann
zeta-function} (2nd edition),  Oxford University Press, Oxford, 1986.


\endRefs
\vskip2cm

\enddocument

\bye